\documentclass[
11pt,
tightenlines,
twoside,
onecolumn,
nofloats,
nobibnotes,
nofootinbib,
superscriptaddress,
noshowpacs,
centertags]
{revtex4-2}
\usepackage{ljm}
\usepackage{amsmath}
\usepackage{amssymb}
\usepackage[ddmmyyyy,hhmmss]{datetime}
\usepackage{enumitem}
\usepackage{url}
\usepackage{blindtext}

\newcommand{\boldm}[1] {\mathversion{bold}#1\mathversion{normal}}

\begin{document}

\titlerunning{Improved MLE of ARMA models}
\authorrunning{L. Di Gangi at al.}

\title{Improved Maximum Likelihood Estimation of ARMA Models}

\author{Leonardo Di Gangi}
\email[E-mail: ]{leonardo.digangi@unifi.it}
\affiliation{Global Optimization Laboratory, DINFO, Universit\`{a} degli Studi di Firenze}

\author{Matteo Lapucci}
%\email[E-mail: ]{matteo.lapucci@unifi.it}
\affiliation{Global Optimization Laboratory, DINFO, Universit\`{a} degli Studi di Firenze}

\author{Fabio Schoen}
%\email[E-mail: ]{fabio.schoen@unifi.it}
\affiliation{Global Optimization Laboratory, DINFO, Universit\`{a} degli Studi di Firenze}

\author{Alessio Sortino}
%\email[E-mail: ]{alessio.sortino@unifi.it}
\affiliation{Global Optimization Laboratory, DINFO, Universit\`{a} degli Studi di Firenze}

%\firstcollaboration{(Submitted by A.~I. Volodin)} % Add if you know submitter.

%\received{November 12, 2021; revised November 20, 2021; accepted
%December 09, 2021}

\begin{abstract}
In this paper we propose a new optimization  model for maximum
likelihood estimation of causal and invertible ARMA models. Through
a set of numerical experiments we show how our proposed model
outperforms, both in terms of quality of the fitted model as well as
in the computational time, the classical estimation procedure based
on Jones reparametrization. We also propose a regularization term in
the model and we show how this addition improves the out of sample
quality of the fitted model. This improvement is achieved thanks to
an increased penalty on models close to the non causality or non
invertibility boundary.
\end{abstract}



Paper accepted for publication in \textit{Lobachevskii Journal of Mathematics}

\maketitle

\section{Introduction}
A zero mean ARMA process of order $(p, q)$ is defined through the following stochastic difference equation:
\begin{equation}
Y_t- \phi_1 Y_{t-1}- \ldots-\phi_p Y_{t-p} =\theta_1 \epsilon_{t-1}+\ldots+\theta_q \epsilon_{t-q}+ \epsilon_t, \quad \epsilon_t \sim \mathcal{WN}(0,\sigma^2),
\label{eq:ARMA equation}
\end{equation}
or in compact form \cite{box2015time} as $
\Phi(B)Y_t=\Theta(B)\epsilon_t, $ where $\Phi(\cdot)$ and
$\Theta(\cdot)$ are the $p$th and $q$th-degree polynomials
\begin{align}
\Phi(z) & = 1- \phi_1 z- \ldots- \phi_p z^p \label{eq: autoregressive polynomial}, \\
\Theta(z) & = 1+ \theta_1 z+ \ldots+ \theta_p z^q, \label{eq: moving-average polynomial}
\end{align}
and $B$ is the backward shift operator (see \cite{box2015time, brockwell1991time,hamilton1994time}). In Equation \eqref{eq:ARMA equation}, $\phi = (\phi_1, \dots, \phi_p)$ are the parameters concerning the autoregressive part, while analogously $\theta = (\theta_1,\dots, \theta_q)$ are the parameters of the moving average component. As it is typically assumed, the error terms in Equation \eqref{eq:ARMA equation} are modeled as a zero mean Gaussian white noise process of variance $\sigma^2$.

The interest towards this class of statistical  models is justified
by their employment in a multitude of fields like business planning,
finance, transportation systems, demography and medicine. With
special reference to real-time forecasting systems, it is very
important to develop computationally efficient estimation methods
focused on improving the numerical stability of the related fitting
procedure and the predictive ability of the ARMA models.

ARMA models estimation has a very long  history
\cite{aigner1971compendium,ansley1979algorithm,box2015time,dent1977computation,gardner1980algorithm,
hannan1982recursive,harvey1979maximum,newbold1974exact}. Maximum
likelihood estimation is usually performed for its advantageous
asymptotic properties. A closed form expression of the ARMA exact
likelihood function was firstly given in \cite{newbold1974exact}.
Afterwards, the focus shifted to finding expressions of the exact
likelihood being more suitable for its computation
\cite{ansley1979algorithm, dent1977computation}. Finally, in the
late 70's, the computational advantages of computing the exact
likelihood by means of Kalman Filter \cite{kalman1960new} have been
pointed out in \cite{harvey1979maximum}. To date, Kalman Filter
algorithm, initialized according to the Gardner method
\cite{gardner1980algorithm}, represents the state-of-the-art of the
methods employed to compute the exact likelihood.

As it is usually required in forecasting  applications, the
estimation of $(\phi, \theta)$ needs to take into account the
causality and invertibility conditions \cite{brockwell1991time}
which act like constraints in the search space. These constraints
are usually handled by means of the Jones reparametrization
\cite{jones1980maximum} which converts the original constrained ARMA
estimation problem into an unconstrained one.

In this paper we propose to fit causal and  invertible ARMA models
by exact maximum likelihood estimation avoiding the employment of
the Jones reparametrization \cite{jones1980maximum}. This is
achievable solving a bound constrained optimization problem. The
benefits of our formulation are both lower computational fitting
times and better numerical stability w.r.t.\ the classical
unconstrained approach. Furthermore, we propose the addition of a
quadratic regularization term to the ARMA exact likelihood function.
This term improves the predictive ability of the fitted ARMA models.

The rest of the paper is organized as follows.  Section 2 contains a
review of the Jones reparametrization method. In Section 3 the
notion of closeness of $(\phi, \theta)$ to the feasibility boundary
is defined. In Section 4 our bound constrained maximum likelihood
estimation approach is provided. In Section 5, extensive
computational experiments which assess the reliability of the
proposed method are reported. Finally, the overall conclusions are
remarked in Section 6.

\section{Jones reparametrization}

When causality and invertibility conditions \cite{brockwell1991time}
hold,  the parameters $\phi = (\phi_1, \dots, \phi_p)$ and $\theta =
(\theta_1,\dots, \theta_q)$ are constrained to belong to the set
$S_p \times S_q$, corresponding to the polynomial operator root
conditions
\begin{align}
S_p &= \{ \phi \in \mathbb{R}^{p}\mid  1- \phi_1 z- \ldots- \phi_p z^p \neq 0 \;\forall\, z\in\mathbb{C} \text{ s.t. } |z|\leq 1  \} \label{eq: feasible autoregressive space} \\
S_q &= \{ \theta \in \mathbb{R}^{q}\mid  1+ \theta_1 z+ \ldots+ \theta_p z^q \neq 0\; \forall\, z\in\mathbb{C} \text{ s.t. } |z|\leq 1 \label{eq: feasible moving-average space}  \rbrace.
\end{align}

These feasible sets are easily identified for $p \leq 2$ and $q \leq
2$,  but for $k > 2$ the form of $S_k$ becomes complicated and for
$k > 4$ the polynomial Equations \eqref{eq: feasible autoregressive
space}, \eqref{eq: feasible moving-average space} cannot be solved
analytically \cite{marriott1995bayesian}. The geometry of the
feasible set $S_p \times S_q$ is described in detail in
\cite{combettes1992best,picinbono1986some,shlien1985geometric}. To
circumvent the problem of dealing with constraints \eqref{eq:
feasible autoregressive space} and \eqref{eq: feasible
moving-average space} Barndorff-Nielsen and Schou
\cite{barndorff1973parametrization} reparametrize $\phi= (\phi_1,
\dots, \phi_p)$ in terms of the partial autocorrelations $\rho =
(\rho_1, \dots, \rho_p)$  by means of the one-to-one continously
differentiable Levinson mapping $\Upsilon(\cdot)$:
\begin{equation}
\label{eq: levinson transformation}
%\begin{aligned}
\phi_k^{(k)}  = \rho_k, \quad k=1,\dots,p,\quad
 \phi_i^{(k)}  =
\phi_i^{(k-1)} - \rho_k \phi_{k-i}^{(k-1)}, \quad i=1,\dots,k-1.
%\end{aligned}
\end{equation}
In \eqref{eq: levinson transformation}, causality is simply obtained
by  $\rho_k \in \; (-1,1) \; \; \forall k = 1, \dots,p$. Jones
\cite{jones1980maximum} introduces an additional mapping $J:
\mathbb{R}^p \to (-1,1)^p$, which allows to formulate the original
problem as an unconstrained optimization problem introducing
variables $u_k, \; k=1, \dots,p$:
\begin{equation}
\label{eq: jones reparametrization}
\rho_k = \frac{1 - \exp{(-u_k)}}{1 + \exp{(-u_k)}}, \quad k = 1, \dots, p.
\end{equation}

Similar transformations can also be employed for the moving average
parameters  $\theta= (\theta_1, \dots, \theta_q)$ in order to
guarantee the invertibility condition. By writing the moving average
polynomial \eqref{eq: moving-average polynomial} for the negative
vector of MA parameters, $-\theta$, we get
\begin{equation}
\Theta(z)  = 1- (-\theta_1) z- \ldots- (-\theta_q) z^q,
\label{eq: moving-average polynomial II}
\end{equation}
and the following can be deduced
\begin{equation}
\label{eq: levinson transformation MA}
\begin{aligned}
\theta_k^{(k)}  &=  b_k, \quad k=1,\dots,q,\\
& \theta_i^{(k)}  = \theta_i^{(k-1)} + b_k \theta_{k-i}^{(k-1)}, \quad i=1,\dots,k-1,
\end{aligned}
\end{equation}
where the variables $b_k \in \; (-1,1) \; \; \forall k = 1, \dots,q$. Jones reparametrization for the moving average part is equivalent to \eqref{eq: jones reparametrization}:
\begin{equation}
\label{eq: jones reparametrization MA}
b_k = \frac{1 - \exp{(-w_k)}}{1 + \exp{(-w_k)}}, \quad k = 1, \dots, q.
\end{equation}
In \cite{jones1980maximum}, the variables $b_k$ are called partial
moving average coefficients. The optimization of the exact
loglikelihood in the causal and invertible feasible space is now
carried out with respect to the variables $u=(u_1,\dots,u_p) \in
\mathbb{R}^p$ and $w=(w_1,\dots,w_q) \in \mathbb{R}^q$.

Note that $\phi=\Upsilon (\rho)$, while $\theta=-\Upsilon(b)$.  In
fact, for any $u$ and $w$, the evaluation of the exact likelihood
function in a causal and invertible feasible point can be computed
by means of the transformations \eqref{eq: levinson transformation},
\eqref{eq: jones reparametrization}, \eqref{eq: levinson
transformation MA}, \eqref{eq: jones reparametrization MA}, and the
Kalman recursions. Inverse Jones transformations are easily found by
solving \eqref{eq: jones reparametrization}, \eqref{eq: jones
reparametrization MA} respectively for $u_k, \; k=1, \dots,p$ and
$w_k, \; k=1, \dots,q$. On the other hand, Monhan
\cite{monahan1984note} derives the expression of the inverse
transformation  $\Upsilon^{-1}(\cdot)$ of \eqref{eq: levinson
transformation} which equivalently can be extended for the moving
average part \eqref{eq: levinson transformation MA}.

\section{Closeness to the Feasiblity Boundary}
\label{sec:boundary} In this Section, the notion of closeness of a
feasible  point $(\phi, \theta) \in S_p \times S_q$ to the set
$\partial S_p \times \partial S_q$, i.e.\  the boundary of the
invertibility and causality regions, is formalized. This will be
useful later in this work, when investigating the relation between
the closeness to the boundary and the numerical stability during the
optimization of the Gaussian ARMA exact log-likelihood function.

It is partially documented \footnote{see, e.g.,
\url{https://www.rdocumentation.org/packages/stats/versions/3.6.2/topics/KalmanLike}
and \url{https://bugs.r-project.org/bugzilla/show_bug.cgi?id=14682}}
that log-likelihood evaluation by Kalman filter may fail when a
point $(\phi, \theta)$ is close to the causality boundary.
Furthermore, it is well known that closeness to the non-invertible
region is problematic due to the presence of the so-called pile-up
effect \cite{kang, pile_up, sargan1983maximum}. Indeed, when the
true parameter of an MA$(1)$ process is close to unity, the model
can be estimated to be non-invertible with a unit root even when the
true process is invertible, with a considerably high probability in
a finite sample. Ansley and Newbold \cite{ansley1980finite} confirm
the presence of such effect in ARMA models too.

Inspired by the method of McLeod and Zhang \cite{mcleod_2006} for
testing the presence of a parameter estimate on the boundary of an
MA$(q)$ model, we define the closeness of a point $(\phi, \theta)$
to the boundary of the invertible and the causal-stationary regions
exploiting the parametrization of an ARMA$(p,q)$ in terms of $\rho$
and $b$:
\begin{align*}
(\phi, \theta)  & = \left (\Upsilon(\rho), -\Upsilon(b) \right ), \\
(\phi, \theta)  & \in S_p \times S_q \iff  (\rho, b) \in (-1, 1)^p \times (-1, 1)^q.
\end{align*}
$\Upsilon(\cdot)$ is not one-to-one on the hypercube boundary
\cite{barndorff1973parametrization}. However,  as elegantly shown in
\cite{mcleod_2006}, $\Upsilon(\cdot)$ maps the boundary of
$(-1,1)^p$ onto $\partial S_p$. Since $\Upsilon(\cdot)$ is a
continuously differentiable function in $[-1,1]^p$, the closeness of
an estimate $\phi \in S_p$ to the non causal-stationary boundary
$\partial S_p$ can be defined respectively in terms of the partial
autocorrelations $\rho$. The same reasoning holds for the moving
average part.

As reported in \cite{mcleod_2006}, $\phi \in \partial S_p$ if and
only if $\|\rho\|_\infty= 1$   and similarly $\theta \in \partial
S_q$ if and only if $\|b\|_\infty= 1$. Now, by fixing a threshold
parameter $\tau > 0$, closeness of  $(\phi, \theta) =
(\Upsilon(\rho), -\Upsilon(b)) \in S_p \times S_q$ to the boundary
$\partial S_p \times \partial S_q$ is defined by the following
conditions:
\begin{enumerate}[label=(\roman*)]
    \item $(\phi, \theta) \in S_p \times S_q$ is close to $\partial S_p$ if and only if $1 - \|\rho\|_\infty < \tau$; \label{eq: closeness to AR boundary}
    \item $(\phi, \theta) \in S_p \times S_q$ is close to $\partial S_q$ if and only if $1-\|b\|_\infty < \tau$; \label{eq: closeness to MA boundary}
    \item $(\phi, \theta) \in S_p \times S_q$ is close to both $\partial S_p$ and $\partial S_q$ if and only if $1 - \|\rho\|_\infty < \tau$ and $1-\|b\|_\infty < \tau$. \label{eq: closeness to AR-MA boundary}
\end{enumerate}
A point $(\phi, \theta) \in S_p \times S_q$ which does not satisfy
any of the above  conditions \ref{eq: closeness to AR boundary},
\ref{eq: closeness to MA boundary}, \ref{eq: closeness to AR-MA
boundary} is defined as a strictly feasible point of $S_p \times
S_q$.

\section{The Proposed Approach}
We propose to fit causal and invertible ARMA$(p,q)$ models by
solving the following bound constrained optimization problem
\begin{equation}
\label{eq: our fitting variant}
\begin{aligned}
\max_{\rho, b, \sigma^2} \;& \ell\left (\Upsilon(\rho),- \Upsilon(b), \sigma^2 \right ) \\
\text{s.t. }& \rho \in \left [-1 + \varepsilon, 1 - \varepsilon \right ]^p,\quad b\in \left [-1 + \varepsilon, 1 - \varepsilon \right ]^q,\quad \sigma\in\mathbb{R}_+ .
\end{aligned}
\end{equation}

Optimizing w.r.t.\ the partial autocorrelation and the partial moving average coefficients avoids the use of the Jones reparametrization \eqref{eq: jones reparametrization}, \eqref{eq: jones reparametrization MA}. Note that this formulation cuts off a small part of the feasible space $S_p \times S_q$. However, as highlighted by thorough numerical experiments that we will describe in the following Section, our formulation provides some nice advantages:
\begin{itemize}
    \item it allows to save a significant amount of running time, as there is no more the need to compute equations \eqref{eq: jones reparametrization} and \eqref{eq: jones reparametrization MA} $p$ and $q$ times respectively, each time the log-likelihood has to be computed during the optimization process (note that every gradient computation by finite differences requires $2(p+q)$ objective evaluations);
    \item it allows to avoid solutions too close to the feasibility boundary that typically lead to numerical errors.
\end{itemize}

We furthermore propose to include in the objective function of
Problem \eqref{eq: our fitting variant} a  Tikhonov regularization
term:
\begin{equation}
\label{eq: regularized fitting variant}
\begin{aligned}
\max_{\rho, b, \sigma^2} \;& \ell \left (\Upsilon(\rho),- \Upsilon(b), \sigma^2 \right ) - \lambda(||\rho||_{2}^2 + ||b||_{2}^2 ) \\
\text{s.t. }& \rho \in \left [-1 + \varepsilon, 1 - \varepsilon \right ]^p,\quad b\in \left [-1 + \varepsilon, 1 - \varepsilon \right ]^q,\quad \sigma\in\mathbb{R}_+ .
\end{aligned}
\end{equation}
We will experimentally show in the following that, in our context,
this term not only discourages solutions close to  the feasibility
boundary, but it also improves the predictive ability of ARMA
models.

\section{Computational Experiments}
In what follows the approximation parameter $\varepsilon$ is set to
$10^{-2}$; we fixed the closeness parameter  $\tau=2\varepsilon$ in
\ref{eq: closeness to AR boundary}, \ref{eq: closeness to MA
boundary}, \ref{eq: closeness to AR-MA boundary}, so that it is
still possible for models \eqref{eq: our fitting variant} and
\eqref{eq: regularized fitting variant} to produce points that are
close to the border of the original feasible set.

All the experiments have been performed on a dataset of
synthetically generated time series.  We simulated a total of 2250
time series of different length $l \in \{100, 1000, 10000\}$ from
ARMA $(p,q)$ Gaussian processes up to a maximum order $(p,q)$ of
$(5,5)$ and standard deviation $\sigma\in \{0.01, 0.1, 1\}$.

Specifically, for a given a combination of length, order and
standard deviation, we generated  10 time series, each one
representing a finite realization of a particular ARMA process with
its structural autoregressive and moving average parameters $(\phi,
\theta)$. Each pair $(\phi, \theta)$ is selected according to the
methodology described in \cite{jones1987randomly}. This methodology
allows to choose $(\phi, \theta)$ from a uniform distribution over
the feasible set $S_p \times S_q$.

Firstly, we are interested in establishing the differences between
solving problem \eqref{eq: our fitting variant}  and the
unconstrained one, based on Jones reparametrization, both from the
standpoints of computational times and numerical stability. To this
aim we carried out a multi-start strategy: for each time series, the
fitting process is repeated 30 times from different randomly chosen
starting points. These starting points are again obtained by uniform
sampling over the feasible region. For a fair comparison, the two
considered methods share the sets of starting points.

Secondly, we investigated the prediction performance of ARMA models
close to the boundary.  As usual, the performance is evaluated on a
test set, after fitting on training data. Our test set for each time
series is given by the last three observations (short term
forecasting scenario). Similarly as above, the process of model
estimation and computation of forecasts is repeated 30 times in a
multi-start fashion. Note that, here, ARMA models have been fitted
only by means of the classical Jones methodology. Indeed, our
interest is to characterize both the forecasting performance of
ARMA models close to the border and how frequently they are obtained
in the standard setting.

Our last experiment assesses the impact of the $\ell_2$
regularization term  in the short term forecasting. For each time
series of our dataset, a single starting point to initialize the
optimization is selected. The fitting procedure is then repeated for
different values of the regularization hyperparameter $\lambda$ in
Equation \ref{eq: regularized fitting variant}.

All the experiments were performed on a machine with Ubuntu Server
20.04 LTS OS,   Intel Xeon E5-2430 v2 @ 2.50GHz CPU and 32GB RAM.

\subsection{Fitting Procedure Runtimes}
Our method provides a significant reduction of the computational
time  required to fit a time series with respect to the
unconstrained fitting method of Jones. The time saving is estimated
to be about $24 \%$ in relative terms.

This result is corroborated by the non parametric Wilcoxon
signed-ranks test  \cite{demsar2006statistical,wilcoxon}. We
considered as fitting time for a time series the average runtime of
successful runs (i.e., with no numerical error) of our multi-start
procedure. Results of the Wilcoxon signed-ranks test are reported in
Tables \ref{tab: Times Wilcoxon two sided} and \ref{tab: Times
Wilcoxon one sided}. These results point out that the median of the
differences of fitting times between the two methods can be assumed
to be positive, i.e., the constrained method has significantly lower
fitting times.
\begin{table}[!h]
    \begin{tabular}{|c|c|}
                \hline
                \textbf{Test statistic} & \textbf{P-value} \\
                \hline
                -34.3807 & $<1\mathrm{e}{-5}$ \\
                \hline
    \end{tabular}
\setcaptionmargin{0mm} 
\onelinecaptionsfalse 
\captionstyle{flushleft} 
\caption{Two sided Wilcoxon signed-rank test. Null hypothesis:  the
median of the differences of the computational times
$t_{\text{Jones}}-$ $t_{\text{our}}$ is zero.} \label{tab: Times
Wilcoxon two sided}
\end{table}

\begin{table}[!h]
    \begin{tabular}{|c|c|}
                \hline
                \textbf{Test statistic} & \textbf{P-value} \\
                \hline
                34.3807 &  $<1\mathrm{e}{-5}$ \\
                \hline
    \end{tabular}
\setcaptionmargin{0mm} 
\onelinecaptionsfalse 
\captionstyle{flushleft} 
\caption{One sided Wilcoxon signed-rank test. Null hypothesis:  the
median of the differences of computational times $t_{\text{Jones}}-$
$t_{\text{our}}$ is negative.} \label{tab: Times Wilcoxon one sided}
\end{table}

\subsection{Numerical Instability}
Our fitting method prevents numerical issues during the optimization
process of the ARMA exact likelihood function, thereby ensuring a
higher level of computational stability.
\begin{table}[h!]
    \begin{center}
        \begin{tabular}{l|c|c}
            \textbf{Method} & \textbf{Arithmetic issues} & \textbf{Kalman Filter errors}\\
            \hline
            Our & 0 & 0.06  \\
            Jones reparametrization & 2.65  & 0.22  \\
        \end{tabular}
\setcaptionmargin{0mm} 
\onelinecaptionsfalse 
\captionstyle{flushleft} 
\caption{Occurrence of numerical instability issues per 1000 runs}
\label{tab: numerical instability}
    \end{center}
\end{table}

The employment of the Jones reparametrization, where exponential
operators are present, leads to a non-negligible probability of
arithmetic issues, which almost always are divisions by zero and in
rare cases overflows. Our method does not suffer at all from these
issues.

The most critical errors, that completely undermine the fitting
process,  come from the Kalman Filter recursions. In general, it is
well known that numerical instability often occurs in Kalman
Filtering \cite{tusell2011kalman}, especially related to the
computation of the state covariance matrix.

Our experiments show that the closeness of a point $(\phi, \theta)$
to the  feasibility boundary is related to numerical instability
within the Kalman Filter recursions. In particular, we observed a
total of 19 \texttt{LinAlgError} errors (15 by the classical method,
4 by using our model \eqref{eq: our fitting variant}) because of the
failed convergence of the SVD numerical computation.

In Tables \ref{tab: errors NO custom} and \ref{tab: errors custom} a
detailed description of these errors is reported. The error may be
due to the evaluation of the log-likelihood in that point or the
computation in the same point of the gradient, since it is
approximated by finite differences.
\begin{table}[h!]
    \begin{center}
        \begin{tabular}{|l|r|r|c|c|c|}
            \hline
            \textbf{Model} & \textbf{Length} & {\boldm $\sigma$}  & \textbf{Starting point} & \textbf{Error point} & \textbf{Ground truth point} \\
            \hline
            ARMA$(2,1)$ & 100 & 0.01 & strictly feasible & (iii) & strictly feasible \\
            ARMA$(2,1)$ & 10000 & 0.01 & strictly feasible & (iii) & (i) \\
            ARMA$(2,1)$ & 10000 & 0.01 & (i) & (ii) & strictly feasible \\
            ARMA$(2,1)$ & 100 & 0.1 & strictly feasible & (iii) & strictly feasible \\
            ARMA$(2,1)$ & 100 & 0.1 & (ii) & (i) & strictly feasible \\
            ARMA$(2,1)$ & 100 & 0.1 & strictly feasible & (i) & strictly feasible \\
            ARMA$(2,1)$ & 1000 & 0.1 & strictly feasible & (iii) & strictly feasible \\
            ARMA$(2,1)$ & 10000 & 0.1 & strictly feasible & (iii) & strictly feasible \\
            ARMA$(2,1)$ & 10000 & 0.1 & strictly feasible & (iii) & strictly feasible \\
            ARMA$(2,1)$ & 100 & 1 & (i) & (iii) & strictly feasible \\
            ARMA$(2,1)$ & 1000 & 1 & strictly feasible & (iii) & strictly feasible \\
            ARMA$(2,1)$ & 10000 & 1 & strictly feasible & (iii) & strictly feasible \\
            ARMA$(2,3)$ & 10000 & 1 & (ii) & (iii) & strictly feasible \\
            ARMA$(3,2)$ & 100 & 0.01 & strictly feasible & (iii) & strictly feasible \\
            ARMA$(5,1)$ & 10000 & 1 & strictly feasible & (i) & strictly feasible \\\hline
        \end{tabular}
    	\setcaptionmargin{0mm} 
    	\onelinecaptionsfalse 
    	\captionstyle{flushleft} 
        \caption{Numerical errors in Kalman filtering when using Jones reparametrization. The first three columns contain information about the ARMA process that generated the tested series and the series itself (orders $p$ and $q$, series length, standard deviation of the white noise generator process). The fourth and fifth columns provide details about the optimization run: the starting point and the point where the error has been generated are characterized in terms of closeness to the feasibility boundary, according to the metrics introduced in Section \ref{sec:boundary}. The sixth column provides the same information associated with the parameters of the model employed to generate the series.}
        \label{tab: errors NO custom}
    \end{center}
\end{table}

\begin{table}[h!]
    \begin{center}
        \begin{tabular}{|l|r|r|r|r|r|r}
            \hline
            \textbf{Model} & \textbf{Length} & {\boldm $\sigma$} & \textbf{Start point} & \textbf{Error point} & \textbf{Ground truth point} \\
            \hline
            ARMA$(4,2)$ & 10000 & 1 & strictly feasible & (iii) & strictly feasible \\
            ARMA$(4,4)$ & 1000 & 0.1 & strictly feasible & (iii) & (ii) \\
            ARMA$(5,5)$ & 100 & 0.1 & strictly feasible & (ii) & strictly feasible \\
            ARMA$(5,5)$ & 1000 & 0.1 & strictly feasible & strictly feasible & strictly feasible  \\
            \hline
        \end{tabular}
        \setcaptionmargin{0mm} 
        \onelinecaptionsfalse 
        \captionstyle{flushleft} 
        \caption{Numerical errors in Kalman filtering when using model \eqref{eq: our fitting variant}. The first three columns contain information about the ARMA process that generated the tested series and the series itself (orders $p$ and $q$, series length, standard deviation of the white noise generator process). The fourth and fifth columns provide details about the optimization run: the starting point and the point where the error has been generated are characterized in terms of closeness to the feasibility boundary, according to the metrics introduced in Section \ref{sec:boundary}. The last column provides the same information associated with the parameters of the model employed to generate the series.}
        \label{tab: errors custom}
    \end{center}
\end{table}

Two patterns are clear from Tables \ref{tab: errors NO custom} and \ref{tab: errors custom}. Firstly, the classical method by Jones fails 4 times more frequently than ours. This means that our reformulation protects from the occurrence of most numerical errors. Secondly, regardless of the type of parametrization employed, it is evident that these numerical errors are related to points close the boundary $\partial S_p \times \partial S_q$ of the feasible set. Furthermore, by observing the first column of both tables, it seems that most errors inside the unconstrained framework happen even when fitting low order models.

\subsection{Forecasting with Almost-Border Models}
As reported above, we employed again a multi-start approach to
assess the predictive performance of close to the border ARMA
models. For our analysis, we picked time series having at least one
strictly feasible solution and at least a solution that meets one of
the conditions \ref{eq: closeness to AR boundary}, \ref{eq:
closeness to MA boundary}, \ref{eq: closeness to AR-MA boundary}. In
doing so, we got a total of 614 time series with such features.

When multiple strictly feasible solutions are available,  we
considered the best one according to the exact log-likelihood value.
The same is done when multiple solutions close to the border are
obtained for a single time series. We then computed multi-step ahead
predictions with the two selected models for each time-series.

Differences in predictive performance of these two distinct ARMA
models  are again investigated by means of the Wilcoxon signed-ranks
test \cite{wilcoxon, demsar2006statistical}. We employed the mean
absolute scaled error (MASE) \cite{hyndman2006another} to measure
the accuracy of forecasts. Indeed, the MASE can be used to compare
forecast methods on a single series and, being scale-free, to
compare forecast accuracy across series \cite{anotherhyndman}.

In our experiments, MASE at a given forecast horizon $h$ is computed as \begin{equation}
\text{MASE}(h)=  \frac{1}{h} \frac{\sum_{t = n+1}^{h} |y_t - \hat{y_t}|}{\frac{1}{n-1} \sum_{t=2}^{n}|y_t - y_{t-1}|}.
\label{eq: mase error}
\end{equation}

We also reported the single absolute scaled errors for each different forecast horizon $h$:
\begin{equation}
\text{ScaledError}(h) = \frac{|y_{n+h}-\hat{y}_{n+h}|}{\frac{1}{n-1} \sum_{t=2}^{n}|y_t - y_{t-1}|}.
\label{eq: absolute scaled error}
\end{equation}

\begin{table}[h!]
    \begin{center}
        \begin{tabular}{|l|c|c|}
            \hline
            \textbf{Error} & \textbf{Test statistic} & \textbf{P-value} \\
            \hline
            MASE$(3)$ & -4.23197 & $2.31\mathrm{e}{-5}$  \\
            \hline
            ScaledError$(1)$ & -1.49874 & 0.13394  \\
            \hline
            ScaledError$(2)$ & -1.67521 & 0.09389 \\
            \hline
            ScaledError$(3)$ & -4.35523 & $1.33\mathrm{e}{-5}$ \\
            \hline
        \end{tabular}
    	\setcaptionmargin{0mm} 
    	\onelinecaptionsfalse 
    	\captionstyle{flushleft} 
        \caption{Results from the two-sided Wilcoxon test at different horizons. Null hypothesis: the median of the differences of the MASE errors, $\text{MASE}_{\text{border}}- \text{MASE}_{\text{strictly feasible}}$, is zero.}
        \label{tab: Two sided Wilcoxon test}
    \end{center}
\end{table}

\begin{table}[h!]
    \begin{center}
        \begin{tabular}{|l|c|c|}
            \hline
            \textbf{Error} & \textbf{Test statistic} & \textbf{P-value} \\
            \hline
            MASE$(3)$ & 4.23197 & $1.16\mathrm{e}{-5}$\\
            \hline
            ScaledError$(1)$ & 1.49874 & 0.06697 \\
            \hline
            ScaledError$(2)$ & 1.67521 & 0.04695 \\
            \hline
            ScaledError$(3)$ & 4.35523 & $<1\mathrm{e}{-5}$ \\
            \hline
        \end{tabular}
    	\setcaptionmargin{0mm} 
    	\onelinecaptionsfalse 
    	\captionstyle{flushleft} 
        \caption{Results from the one-sided Wilcoxon test at different horizons. Null
         hypothesis: the median of the differences of the MASE errors, $\text{MASE}_{\text{border}}- \text{MASE}_{\text{strictly feasible}}$, is negative.}
        \label{tab: One sided Wilcoxon test}
    \end{center}
\end{table}

Results are reported in Tables \ref{tab: Two sided Wilcoxon test}
and \ref{tab: One sided Wilcoxon test}. The observed P-value in the
last row of Table \ref{tab: Two sided Wilcoxon test} evidences that
significant differences exist in forecast accuracy between strictly
feasible ARMA $(p, q)$ models and close-to-the-border ARMA $(p, q)$
models. The significant differences involve only the MASE $(3)$
error and the absolute scaled error at horizon $h=3$: in both cases
the associated P-values are strictly lower than the default
significance level $\alpha=0.05$. Furthermore, for these two metrics
the one-sided test confirms that ARMA models close to the
feasibility boundary perform poorer in terms of the predictive
ability than the strictly feasible ARMA models.

Considering instead the remaining error metrics, results  in Table
\ref{tab: Two sided Wilcoxon test} indicate that at forecast horizon
$h=1$ non substantial difference exists in forecast accuracy between
the two types of ARMA models. Differences in predictive ability
become more evident as the forecast horizon grows. From Table
\ref{tab: Two sided Wilcoxon test} we observe that at horizon 2,
only assuming a significance level $\alpha = 0.1$, it is possible to
deduce a statistically significant difference between the two ARMA
models in forecasting performances.

The main conclusion of this experiment is that ARMA models  close to
the feasibility boundary perform poorer in terms of the predictive
ability than the strictly feasible ARMA models. The practical
meaning of this result is that caution is needed with close to the
border ARMA models when forecasting is required. This is one of the
motivations to modify our fitting model \eqref{eq: our fitting
variant} by adding to the objective an $\ell_2$ penalty term as in
\eqref{eq: regularized fitting variant}. We will discuss in depth
the effects of this modification in the next section.

\subsection{Forecasting with Regularized ARMA models}
The next and final experiment investigates the effect of the
addition of an $\ell_2$-regularization term from a forecasting
accuracy perspective. Different values of the regularization
hyperparameter $\lambda$ in Equation \eqref{eq: regularized fitting
variant} give rise to different ARMA$(p,q)$ models with diverse
forecasting performances.

ARMA models are, in practice, fitted by iterative optimization
algorithms that start at preliminary estimates obtained, for
example, with the well-known Hannan and Rissanen (HR) method
\cite{hannan1982recursive}. We consider this setting to carry out
the experiment, in order to assess the impact of the regularization
term in the common use cases.

The classical Jones fitting method is compared with  models
\eqref{eq: our fitting variant} and \eqref{eq: regularized fitting
variant}, varying the values of the regularization parameter
$\lambda$. For each time series, all optimization algorithms are
started at the same initial point, identified using HR procedure.

We employed the Friedman test
\cite{demsar2006statistical,friedman_1,friedman_2} to  catch the
differences between the methods. The test ranks the fitting methods
for each time series separately, the best performing method (lowest
error) getting the rank of 1, the second best rank 2 and so on. The
null-hypothesis, states that all the fitting methods are equivalent
and so their ranks should be equal. Table \ref{1}
%tab: Average of
%ranks from Friedman test}
reports the average of ranks over all the
time series in our dataset, w.r.t.\ the metrics of interest
\eqref{eq: mase error} and \eqref{eq: absolute scaled error}.

We observe from Table \ref{1} %tab: Average of ranks from Friedman
%test}
that  for the MASE$(3)$ and the absolute scaled error at
horizon $h=3$ the averages of ranks go down until a value of the
hyperparameter $\lambda = 8$. For the other two errors the trend of
the averages of the ranks seems quite stationary: this pattern finds
confirmation from the results of Friedman test as it is shown in
Table \ref{tab: Friedman test}.

\begin{table}[!h]
    \begin{center}
        {
            \begin{tabular}{|l|c|c|c|c|c|c|c|}
                \hline
                \textbf{Error} & \textbf{Jones} & $\boldsymbol{\lambda=0}$ & $\boldsymbol{\lambda=1}$ & $\boldsymbol{\lambda=2}$ & $\boldsymbol{\lambda=4}$ & $\boldsymbol{\lambda=8}$ & $\boldsymbol{\lambda=16}$ \\
                \hline
                MASE$(3)$ &        4.228 &        4.201 &       4.056 &     3.947 &     3.882 &     \textbf{3.825} &     3.862 \\
                \hline
                ScaledError$(1)$ &        4.022 &        3.996 &       4.018 &     3.999 &     3.972 &     \textbf{3.968} &     4.025 \\
                \hline
                ScaledError$(2)$ &        4.082 &        4.095 &       4.01 &     3.972 &     3.958 &     \textbf{3.935} &     3.948 \\
                \hline
                ScaledError$(3)$ &        4.220 &        4.226 &       4.081 &     3.980 &     3.885 &    \textbf{3.798} &     3.809 \\
                \hline
            \end{tabular}
        }
        \setcaptionmargin{0mm} 
        \onelinecaptionsfalse 
        \captionstyle{flushleft} 
        \caption{Average of ranks between different ARMA models performance w.r.t.\ different error metrics.}
        \label{1}%tab: Average of ranks from Friedman test}
    \end{center}
\end{table}

\begin{table}[!h]
    \begin{center}
        \begin{tabular}{|l|c|c|}
            \hline
            \textbf{Error} & \textbf{Test statistic} & \textbf{P-value}  \\
            \hline
            MASE$(3)$ &        78.06724 &     $<1\mathrm{e}{-5}$  \\
            \hline
            ScaledError$(1)$   &        1.57091  & 0.95465  \\
            \hline
            ScaledError$(2)$   &        12.13886 &     0.05894  \\
            \hline
            ScaledError$(3)$   &        94.93939 & $<1\mathrm{e}{-5}$ \\
            \hline
        \end{tabular}
    	\setcaptionmargin{0mm} 
    	\onelinecaptionsfalse 
    	\captionstyle{flushleft} 
        \caption{Results of Friedman test for the difference in forecasting performance of various ARMA models w.r.t.\ different error metrics.}
        \label{tab: Friedman test}
    \end{center}
\end{table}

Friedman test, whose results are reported in Table  \ref{tab:
Friedman test}, suggests that the forecasting performance of the
considered fitting models statistically differ (assuming a
significance level of $\alpha= 0.1$) for all the errors except for
the absolute scaled forecasting error at horizon $h=1$.

Therefore, based on these results we considered necessary to conduct
post hoc-analysis w.r.t.\ the MASE$(3)$, the absolute scaled
forecasting error at horizon $h=3$ and $h=2$ (although the P-value
in the latter case is not negligible).

Post-hoc analysis is performed by means of the Nemenyi test
\cite{demsar2006statistical,nemenyi1962distribution}. Critical
differences between two generic methods are assessed in terms of the
differences between the averages of the ranks. Results of the
Nemenyi test are reported in Tables \ref{tab: nemenyi MASE},
\ref{tab: nemenyi scaled absolute 2} and \ref{tab: nemenyi scaled
absolute 3}.

\begin{table}[!h]
    \begin{center}
        {
            \begin{tabular}{|l|c|c|c|c|c|c|c|}
                \hline
                {} &  \textbf{Jones} & $\boldsymbol{\lambda=0}$ & $\boldsymbol{\lambda=1}$ & $\boldsymbol{\lambda=2}$ & $\boldsymbol{\lambda=4}$ & $\boldsymbol{\lambda=8}$ & $\boldsymbol{\lambda=16}$ \\
                \hline
                \textbf{Jones} &                   1.00000 &                 0.90000 &                 0.10395 &                 0.00100 &                 0.00100 &                 0.00100 &                  0.00100 \\
                \hline
                $\boldsymbol{\lambda=0}$   &                   0.90000 &                 1.00000 &                 0.26546 &                 0.00154 &                 0.00100 &                 0.00100 &                  0.00100 \\
                \hline
                $\boldsymbol{\lambda=1}$   &                   0.10395 &                 0.26546 &                 1.00000 &                 0.60537 &                 0.10031 &                 0.00630 &                  0.04196 \\
                \hline
                $\boldsymbol{\lambda=2}$   &                   0.00100 &                 0.00154 &                 0.60537 &                 1.00000 &                 0.90000 &                 0.48698 &                  0.82448 \\
                \hline
                $\boldsymbol{\lambda=4}$   &                   0.00100 &                 0.00100 &                 0.10031 &                 0.90000 &                 1.00000 &                 0.90000 &                  0.90000 \\
                \hline
                $\boldsymbol{\lambda=8}$   &                   0.00100 &                 0.00100 &                 0.00630 &                 0.48698 &                 0.90000 &                 1.00000 &                  0.90000 \\
                \hline
                $\boldsymbol{\lambda=16}$  &                   0.00100 &                 0.00100 &                 0.04196 &                 0.82448 &                 0.90000 &                 0.90000 &                  1.00000 \\
                \hline
            \end{tabular}
        }
        \setcaptionmargin{0mm} 
        \onelinecaptionsfalse 
        \captionstyle{flushleft} 
        \caption{Posthoc analysis of the performance forecasting: pairwise comparison of the MASE$(3)$ error.}
        \label{tab: nemenyi MASE}
    \end{center}
\end{table}

\begin{table}[!h]
    \begin{center}
        {
            \begin{tabular}{|l|c|c|c|c|c|c|c|}
                \hline
                {} &  \textbf{Jones} & $\boldsymbol{\lambda=0}$ & $\boldsymbol{\lambda=1}$ & $\boldsymbol{\lambda=2}$ & $\boldsymbol{\lambda=4}$ & $\boldsymbol{\lambda=8}$ & $\boldsymbol{\lambda=16}$ \\
                \hline
                \textbf{Jones} &                   1.00000 &                 0.90000 &                 0.90000 &                 0.60131 &                 0.46951 &                 0.25145 &                  0.37172 \\
                \hline
                $\boldsymbol{\lambda=0}$   &                   0.90000 &                 1.00000 &                 0.82448 &                 0.48264 &                 0.34176 &                 0.16502 &                  0.25839 \\
                \hline
                $\boldsymbol{\lambda=1}$   &                   0.90000 &                 0.82448 &                 1.00000 &                 0.90000 &                 0.90000 &                 0.90000 &                  0.90000 \\
                \hline
                $\boldsymbol{\lambda=2}$   &                   0.60131 &                 0.48264 &                 0.90000 &                 1.00000 &                 0.90000 &                 0.90000 &                  0.90000 \\
                \hline
                $\boldsymbol{\lambda=4}$   &                   0.46951 &                 0.34176 &                 0.90000 &                 0.90000 &                 1.00000 &                 0.90000 &                  0.90000 \\
                \hline
                $\boldsymbol{\lambda=8}$   &                   0.25145 &                 0.16502 &                 0.90000 &                 0.90000 &                 0.90000 &                 1.00000 &                  0.90000 \\
                \hline
                $\boldsymbol{\lambda=16}$  &                   0.37172 &                 0.25839 &                 0.90000 &                 0.90000 &                 0.90000 &                 0.90000 &                  1.00000 \\
                \hline
            \end{tabular}
        }
        \setcaptionmargin{0mm} 
        \onelinecaptionsfalse 
        \captionstyle{flushleft} 
        \caption{Posthoc analysis of the performance forecasting: pairwise comparison of the absolute scaled error at horizon $h =2$.}
        \label{tab: nemenyi scaled absolute 2}
    \end{center}
\end{table}

\begin{table}[!h]
    \begin{center}{
            \begin{tabular}{|l|c|c|c|c|c|c|c|}
                \hline
                {} &  \textbf{Jones} & $\boldsymbol{\lambda=0}$ & $\boldsymbol{\lambda=1}$ & $\boldsymbol{\lambda=2}$ & $\boldsymbol{\lambda=4}$ & $\boldsymbol{\lambda=8}$ & $\boldsymbol{\lambda=16}$ \\
                \hline
                \textbf{Jones} &                   1.00000 &                 0.90000 &                 0.31753 &                 0.00357 &                 0.00100 &                 0.00100 &                  0.00100 \\
                \hline
                $\boldsymbol{\lambda=0}$   &                   0.90000 &                 1.00000 &                 0.27263 &                 0.00259 &                 0.00100 &                 0.00100 &                  0.00100 \\
                \hline
                $\boldsymbol{\lambda=1}$   &                   0.31753 &                 0.27263 &                 1.00000 &                 0.67435 &                 0.03709 &                 0.00100 &                  0.00100 \\
                \hline
                $\boldsymbol{\lambda=2}$   &                   0.00357 &                 0.00259 &                 0.67435 &                 1.00000 &                 0.73116 &                 0.07136 &                  0.11154 \\
                \hline
                $\boldsymbol{\lambda=4}$   &                   0.00100 &                 0.00100 &                 0.03709 &                 0.73116 &                 1.00000 &                 0.80825 &                  0.90000 \\
                \hline
                $\boldsymbol{\lambda=8}$   &                   0.00100 &                 0.00100 &                 0.00100 &                 0.07136 &                 0.80825 &                 1.00000 &                  0.90000 \\
                \hline
                $\boldsymbol{\lambda=16}$  &                   0.00100 &                 0.00100 &                 0.00100 &                 0.11154 &                 0.90000 &                 0.90000 &                  1.00000 \\
                \hline
            \end{tabular}
        }
        \setcaptionmargin{0mm} 
        \onelinecaptionsfalse 
        \captionstyle{flushleft} 
        \caption{Posthoc analysis of the performance forecasting: pairwise comparison of the absolute scaled error at horizon $h =3$.}
        \label{tab: nemenyi scaled absolute 3}
    \end{center}
\end{table}

Regarding the absolute scaled error at horizon $h=2$, results from
the Nemenyi  test indicate no significant differences between the
fitting methods in terms of the forecasting performances. All the
P-values reported in Table \ref{tab: nemenyi scaled absolute 2} are
greater than 0.1.

On the other end, results about absolute scaled error at horizon
$h=3$  and the MASE$(3)$ are equivalent. By observing both Table
\ref{tab: nemenyi MASE} and Table \ref{tab: nemenyi scaled absolute
3}, no significant difference is found between the two non
regularized methods. Furthermore, no significant differences in
forecasting performance have been identified between both the non
regularized methods and the regularized one with $\lambda=1$.

Instead, stronger regularization leads to significantly better
forecasts  w.r.t.\ the non  regularized methods. Forecasting
performance, as mentioned above, starts to deteriorate as the
regularization hyperparameter grows to $\lambda=16$. In summary, the
constrained fitting method with regularization leads to causal and
invertible ARMA models with better short term predictive ability
than the non regularized ones.

\section{Conclusions}
Fitting causal and invertible ARMA models by constrained
optimization in the partial  autocorrelation and partial
moving-average coefficients space has several advantages w.r.t.\ the
classical unconstrained approach based on the Jones
reparametrization. First of all, we observed that our approach leads
to a significant reduction of the fitting times. Moreover,
almost-border solutions are often avoided. Such solutions, as
further experiments highlight, are bad both because they lead to
numerical errors during the optimization of the ARMA exact
log-likelihood and because they do not perform well at forecasting.

Based on these results we proposed $\ell_2$-regularization to
discourage  almost-border solutions. As non parametric statistical
tests assess, $\ell_2$-regularization also improves the short term
forecasting performances of causal and invertible ARMA models.


\begin{thebibliography}{99}
	\bibitem{aigner1971compendium}
	\refitem{article}
	D.~J. Aigner,
	\textquotedblleft{A compendium on estimation of the autoregressive moving average model from the series data}\textquotedblright, International Economic Review, pages 348--371, 1971.
	
	\bibitem{ansley1979algorithm}
	\refitem{article}
	C.~F. Ansley, \textquotedblleft{An algorithm for the exact likelihood of a mixed
		autoregressive-moving average process}\textquotedblright, Biometrika, 66(1):59--65, 1979.
	
	\bibitem{ansley1980finite}
	\refitem{article}
	C.~F. Ansley and P.~Newbold,
	\textquotedblleft{Finite sample properties of estimators for autoregressive moving average models}\textquotedblright,
	Journal of Econometrics, 13(2):159--183, 1980.
	
	\bibitem{barndorff1973parametrization}
	\refitem{article}
	O.~Barndorff-Nielsen and G.~Schou,
	\textquotedblleft{On the parametrization of autoregressive models by partial autocorrelations}\textquotedblright, 
	Journal of multivariate Analysis, 3(4):408--419, 1973.
	
	\bibitem{box2015time}
	\refitem{book}
	G.~E. Box, G.~M. Jenkins, G.~C. Reinsel, and G.~M. Ljung, \emph{Time series analysis: forecasting and control} (John Wiley \& Sons, 2015).
	
	\bibitem{brockwell1991time}
	\refitem{book}
	P.~J. Brockwell, R.~A. Davis, and S.~E. Fienberg, \emph{Time series: theory and methods: theory and methods} (Springer Science \& Business Media, 1991).
	
	\bibitem{combettes1992best}
	\refitem{article}
	P.~L. Combettes and H.~J. Trussell,\textquotedblleft{Best stable and invertible approximations for ARMA systems}\textquotedblright, IEEE Transactions on signal processing, 40(12):3066--3069, 1992.
	
	\bibitem{demsar2006statistical}
	\refitem{article}
	J.~Dem{\v{s}}ar,
	\textquotedblleft{Statistical comparisons of classifiers over multiple data sets}\textquotedblright, The Journal of Machine Learning Research, 7:1--30, 2006.
	
	\bibitem{dent1977computation}
	\refitem{article}
	W.~Dent,
	\textquotedblleft{Computation of the exact likelihood function of an arima process}\textquotedblright, Journal of Statistical Computation and Simulation, 5(3):193--206, 1977.
	
	\bibitem{friedman_1}
	\refitem{article}
	M.~Friedman,
	\textquotedblleft{The use of ranks to avoid the assumption of normality implicit in the analysis of variance}\textquotedblright, Journal of the american statistical association, 32(200):675--701, 1937.
	
	\bibitem{friedman_2}
	\refitem{article}
	M.~Friedman,
	\textquotedblleft{A comparison of alternative tests of significance for the problem of m ranking}\textquotedblright, The Annals of Mathematical Statistics, 11(1):86--92, 1940.
	
	\bibitem{gardner1980algorithm}
	\refitem{article}
	G.~Gardner, A.~C. Harvey, and G.~D. Phillips,
	\textquotedblleft{Algorithm as 154: An algorithm for exact maximum likelihood estimation of autoregressive-moving average models by means of kalman filtering}\textquotedblright, Journal of the Royal Statistical Society. Series C (Applied Statistics), 29(3):311--322, 1980.
	
	\bibitem{hamilton1994time}
	\refitem{book}
	J.~D. Hamilton, \emph{Time series analysis, volume~2} (Princeton New Jersey, 1994).
	
	\bibitem{hannan1982recursive}
	\refitem{article}
	E.~J. Hannan and J.~Rissanen, \textquotedblleft{Recursive estimation of mixed autoregressive-moving average order}\textquotedblright, Biometrika, 69(1):81--94, 1982.
	
	\bibitem{harvey1979maximum}
	\refitem{article}
	A.~C. Harvey and G.~D. Phillips, \textquotedblleft{Maximum likelihood estimation of regression models with autoregressive-moving average disturbances}\textquotedblright, Biometrika, 66(1):49--58, 1979.
	
	\bibitem{anotherhyndman}
	\refitem{article}
	R.~J. Hyndman et~al, \textquotedblleft{Another look at forecast-accuracy metrics for intermittent demand}\textquotedblright, Foresight: The International Journal of Applied Forecasting, 4(4):43--46, 2006.
	
	\bibitem{hyndman2006another}
	\refitem{article}
	R.~J. Hyndman and A.~B. Koehler, \textquotedblleft{Another look at measures of forecast accuracy}\textquotedblright, Another look at measures of forecast accuracy, International journal of forecasting, 22(4):679--688, 2006.
	
	\bibitem{jones1987randomly}
	\refitem{article}
	M.~Jones, \textquotedblleft{Randomly choosing parameters from the stationarity and invertibility region of autoregressive--moving average models}\textquotedblright, Journal of the Royal Statistical Society: Series C (Applied Statistics), 36(2):134--138, 1987.
	
	\bibitem{jones1980maximum}
	\refitem{article}
	R.~H. Jones, \textquotedblleft{Maximum likelihood fitting of arma models to time series with missing observations}\textquotedblright, Technometrics, 22(3):389--395, 1980.
	
	\bibitem{kalman1960new}
	\refitem{article}
	R.~E. Kalman, \textquotedblleft{A new approach to linear filtering and prediction problems}\textquotedblright, 1960.
	
	\bibitem{kang}
	\refitem{article}
	K.~M. Kang, \textquotedblleft{A comparison of estimators for moving average processes}\textquotedblright, Unpublished Paper, Australian Bureau of Statistics, 1975.
	
	\bibitem{pile_up}
	\refitem{article}
	C.-J. Kim and J.~Kim, \textquotedblleft{Thepile-up problem'in trend-cycle decomposition of real gdp: Classical and bayesian perspectives}\textquotedblright, 2013.
	
	\bibitem{marriott1995bayesian}
	\refitem{article}
	J.~Marriott \textquotedblleft{Bayesian analysis of arma processes: Complete sampling-based inferences under full likelihood}\textquotedblright, Bayesian Statistics and Econometrics: Essays in Honor of Arnold Zellner, 1995.
	
	\bibitem{monahan1984note}
	\refitem{article}
	J.~F. Monahan, \textquotedblleft{A note on enforcing stationarity in autoregressive-moving average models}\textquotedblright, Biometrika, 71(2):403--404, 1984.
	
	\bibitem{nemenyi1962distribution}
	\refitem{article}
	P.~Nemenyi,
	\textquotedblleft{Distribution-free multiple comparisons}\textquotedblright, In Biometrics, volume~18, page 263. International Biometric Soc 1441 I ST, NW, SUITE 700, WASHINGTON, DC 20005-2210, 1962.
	
	\bibitem{newbold1974exact}
	\refitem{article}
	P.~Newbold,
	\textquotedblleft{The exact likelihood function for a mixed autoregressive-moving average process}\textquotedblright, Biometrika, 61(3):423--426, 1974.
	
	\bibitem{picinbono1986some}
	\refitem{article}
	B.~Picinbono and M.~Benidir,
	\textquotedblleft{Some properties of lattice autoregressive filters}\textquotedblright, IEEE transactions on acoustics, speech, and signal processing, 34(2):342--349, 1986.
	
	\bibitem{sargan1983maximum}
	\refitem{article}
	J.~D. Sargan and A.~Bhargava,
	\textquotedblleft{Maximum likelihood estimation of regression models with first order moving average errors when the root lies on the unit circle}\textquotedblright, Econometrica: Journal of the Econometric Society, pages 799--820, 1983.
	
	\bibitem{shlien1985geometric}
	\refitem{article}
	S.~Shlien,
	\textquotedblleft{A geometric description of stable linear predictive coding digital filters (corresp.)}\textquotedblright, IEEE Transactions on information theory, 31(4):545--548, 1985.
	
	\bibitem{tusell2011kalman}
	\refitem{article}
	F.~Tusell et~al,
	\textquotedblleft{Kalman filtering in r}\textquotedblright, Journal of Statistical Software, 39(2):1--27, 2011.
	
	\bibitem{wilcoxon}
	\refitem{article}
	F.~Wilcoxon,
	\textquotedblleft{Individual comparisons by ranking methods}\textquotedblright, In Breakthroughs in statistics, pages 196--202. Springer, 1992.
	
	\bibitem{mcleod_2006}
	\refitem{article}
	Y.~Zhang and A.~I. McLeod,
	\textquotedblleft{Fitting ma (q) models in the closed invertible region}\textquotedblright, Statistics \& probability letters, 76(13):1331--1334, 2006.
	\thispagestyle{empty}
\end{thebibliography}
\end{document}